\documentclass[11pt,reqno,oneside]{amsart}
 \usepackage{latexsym}
 \usepackage{amssymb}
 \usepackage{amsfonts}
\usepackage{graphics}
\usepackage{mathcomp}
 \author[Mathonet]{P. Mathonet}
\thanks{Institute of mathematics, Grande Traverse, 12 - B37, B-4000
 Li\`ege, Belgium\\\indent email : P.Mathonet@ulg.ac.be, F.Radoux@ulg.ac.be}
\author[Radoux]{F. Radoux}
\date{\today} 
\title[Natural and Projectively Equivariant Quantizations]{Natural and Projectively Equivariant Quantizations by means of Cartan Connections}

\newtheorem{lem}{Lemma}
\newtheorem{thm}[lem]{Theorem}
\newtheorem{prop}[lem]{Proposition}

\theoremstyle{remark}

\theoremstyle{definition}
\newtheorem{defi}{Definition}

\newcommand{\R}{\mathbb{R}}

\newcommand{\N}{\mathbb{N}}

\newcommand{\g}{\mathfrak{g}}

\newcommand{\h}{\mathfrak{h}}

\begin{document}
\begin{abstract}
The existence of a natural and projectively equivariant quantization in the
sense of Lecomte \cite{Leconj} was proved recently by M. Bordemann \cite{Bor},
using the framework of Thomas-Whitehead connections. We give a new proof of
existence using the notion of Cartan projective connections and we obtain
an explicit formula in terms of these connections. Our method yields the
existence of a projectively equivariant quantization if and only if an
$sl(m+1,\R)$-equivariant quantization exists in the flat situation in the sense of \cite{LO}, thus solving one of
the problems left open by M. Bordemann.
\end{abstract}
\maketitle
\section{Introduction}
Among the different meanings that the word quantization can assume, one (in
the framework of geometric quantization) is to
think of a quantization procedure as a linear bijection from the space of
classical observables to a space of differential operators acting on wave
functions (see \cite{Woo}).
More precisely, in our setting,  the space of observables (also called the
space of \emph{Symbols}) is the space of smooth functions on
the cotangent bundle $T^*M$ of a manifold $M$, that are polynomial along the
fibres. The space of differential operators $\mathcal{D}_{\frac{1}{2}}(M)$ is
made of differential operators acting on half-densities. It is known that
there is no natural quantization procedure. In other words, the spaces of
symbols and of differential operators are not isomorphic as representations of
$\mathrm{Diff}(M)$.

However, when there is a Lie group $G$ acting on $M$ by local diffeomorphisms,
the action can be lifted to symbols and differential operators and 
one can raise the question of knowing whether these spaces are isomorphic
representations of $G$ or not. This leads to the concept of $G$-equivariant
quantization introduced by P. Lecomte and V. Ovsienko in \cite{LO} : a
$G$-equivariant quantization is a linear bijection from the space of symbols
to the space of differential operators that exchanges the actions of $G$ on
these spaces.

In \cite{LO}, the authors considered the case of the projective group
$PGL(m+1,\R)$ acting on the manifold $M=\R^m$ by linear fractional
transformations. This leads to the notion of projectively equivariant
quantization or its infinitesimal counterpart, 
the $sl(m+1,\R)$-equivariant quantization.
One of the main results in this case is the existence of a projectively
equivariant quantization and its uniqueness, up to some natural normalization
condition.
The authors also showed that their results could be directly generalized to
the case of a manifold endowed with a flat projective structure.    

In \cite{DLO}, the authors considered the group $SO(p+1,q+1)$ acting on the
space $\R^{p+q}$ or on a manifold endowed with a flat conformal structure.
 They extended the problem by considering the space $\mathcal{D}_{\lambda,\mu}$
of differential operators mapping $\lambda$-densities into $\mu$-densities and
a suitable space of symbols $\mathcal{S}_{\mu-\lambda}$.
There again, the result was the existence and uniqueness of a conformally
equivariant quantization provided the shift value $\delta = \mu - \lambda$
does not belong to a set of critical values.
Similar results for other Lie groups $G$ acting on vector spaces or other
types of differential operators  were obtained in \cite{BM,DO,BHMP}.

At that point, all these results were dealing with a manifold endowed with a
 flat structure. It was remarked in \cite{Bou1,Bou2} that the formula for the
 projectively equivariant quantization for differential operators of order two and three
 could be generalized to an arbitrary manifold.
 In these papers, S. Bouarroudj showed how to define a
 quantization map from the space of symbols to
 the space of differential operators, using a torsion-free connection, in
 such a way that the quantization map depends only on the \emph{projective
 class} of the connection (recall that two torsion-free linear connections are
 projectively equivalent if they define the same paths, that is, the same
 geodesics up to parametrization). 

In \cite{Leconj}, P. Lecomte conjectured the existence of a quantization 
procedure depending on a torsion-free connection, that would be
 natural (in all arguments) and that would be left invariant by a projective
 change of connection.

The existence of such a quantization procedure was proved by M. Bordemann 
in \cite{Bor}. 
In order to prove the existence, M. Bordemann used a construction that can be
roughly summarized as follows :  first he associated to each projective class
$[\nabla]$  of torsion-free linear connections on $M$ a unique
 linear connection $\tilde{\nabla}$ on a
 principal line bundle $\tilde{M}\to M$, then he showed how to lift the symbols
 to a suitable space of tensors on $\tilde{M}$, and he eventually applied
 the so-called \emph{Standard ordering}. This construction was later adapted
 in \cite{Hansoul} in order to deal with differential operators acting on
 forms.

The study of the projective equivalence of connections goes back to the
 1920's. At that time, there were two main approaches to the so-called
 \emph{Geometry of paths}.
 The connection used by M. Bordemann is inpired by the approach due
 to T.Y. Thomas \cite{Thomas1}, J.H.C. Whitehead
 \cite{Whitehead} and  O.Veblen \cite{VT}(see also
 \cite{Roberts2,Roberts1,Roberts3} for a modern formulation).

 The second approach, due to
 E. Cartan \cite{Cartan1}, leads to the concept of Cartan projective 
connection, developed in a modern setting by S. Kobayashi and T. Nagano in
 \cite{Kobaproj,Kobabook}.

In this paper, we analyse the existence of a natural and projectively
equivariant quantization map from the space of symbols
$\mathcal{S}_{\mu-\lambda}$ to the space $\mathcal{D}_{\lambda,\mu}$, using
 Cartan connections. We obtain an explicit formula for the quantization map in
 terms of the normal Cartan connection associated to a projective equivalence
 class of torsion free-linear connections. This formula generalizes the one
 given by M. Bordemann in \cite{Bor} and is nothing but the formula for the
 flat case given in \cite{DO}.
In particular, we show that the natural and projectively equivariant
 quantization map exists if and only if an $sl(m+1)$-equivariant quantization
 exists in the flat case, thus solving a problem left open by M. Bordemann.

We believe that our methods will apply in order to solve the problem in the
conformal situation or in order to define a projectively equivariant symbol
calculus for other types of differential operators.
\section{Problem setting}
For the sake of completeness, we briefly recall in this section the
definitions of tensor densities, differential operators and symbols. Then we
set the problem of existence of projectively equivariant natural quantizations.
Throughout this note, we denote by $M$ a smooth, Hausdorff and
second countable manifold of dimension $m$.
\subsection{Tensor densities}
The vector bundle of tensor densities $F_{\lambda}(M)\to M$ is a line
bundle associated to the linear frame bundle :
\[ F_{\lambda}(M) = P^1M\times_{\rho}\Delta^{\lambda}(R^m),\]
where the representation $\rho$ of the group $GL(m,\R)$ on the one-dimensional
vector space $\Delta^{\lambda}(R^m)$ is given by
\[\rho(A) e = \vert det A\vert^{-\lambda} e,\quad\forall A\in
GL(m,\R),\;\forall e\in \Delta^{\lambda}(\R^m).\]
As usual, we denote by $\mathcal{F}_{\lambda}(M)$ the space of smooth
sections of this bundle. This is the space $C^{\infty}(P^1M,
\Delta^{\lambda}(R^m))_{GL(m,\R)}$ of functions $f$ such that
\[f(u A) = \rho (A^{-1}) f(u)\quad \forall u \in P^1M,\;\forall A\in
GL(m,\R).\]
Since $ F_{\lambda}(M)\to M$ is  associated to $P^1M$, there are natural
actions of $\mathrm{Diff}(M)$ and of $\mathrm{Vect}(M)$ on
$\mathcal{F}_{\lambda}(M)$. For more details, we refer the reader to \cite{DLO,LO}.
\subsection{Differential operators and symbols}
As in \cite{LO,Bor}, we denote by $\mathcal{D}_{\lambda,\mu}(M)$ the space of
differential operators from $\mathcal{F}_{\lambda}(M)$ to
$\mathcal{F}_{\mu}(M)$. The actions of $\mathrm{Vect}(M)$ and
$\mathrm{Diff}(M)$ are  induced by the actions on tensor densities~:
One has
\[(\phi\cdot D)(f) = \phi\cdot(D(\phi^{-1}\cdot f)),\quad\forall f\in
\mathcal{F}_{\lambda}(M), D\in \mathcal{D}_{\lambda,\mu},\mbox{and}\, \phi
\in \mathrm{Diff}(M).\]
The space $\mathcal{D}_{\lambda,\mu}$ is filered by the order of
differential operators. We denote by $\mathcal{D}^k_{\lambda,\mu}$ the
space of differential operators of order at most $k$. It is well-known
that this filtration is preserved by the action of local diffeomorphisms.
The space of \emph{symbols} is then the associated graded space of
$\mathcal{D}_{\lambda,\mu}$.

We denote by $S^l_{\delta}(\R^m)$ the vector space
 $S^l\R^m\otimes\Delta^{\delta}(\R^m)$. There is a natural representation
 $\rho$ of $GL(m,\R)$ on this space (the representation of $GL(m,\R)$ on
 symmetric tensors is the natural one). We then denote by 
$S^l_{\delta}(M)\to M$
 the
vector bundle
\[P^1M\times_{\rho}S^l_{\delta}(\R^m)\to M,\]
and by $\mathcal{S}^l_{\delta}(M)$ the space of smooth sections of
$S^l_{\delta}(M)\to M$, that is, the space $C^{\infty}(P^1M,
S^l_{\delta}(\R^m))_{ GL(m,\R)}$.

 Then if $\delta = \mu
- \lambda$ the \emph{principal symbol operator} $\sigma :
\mathcal{D}^l_{\lambda,\mu}(M)\to \mathcal{S}^l_{\delta}(M)$ commutes with
the action of diffeomorphisms and is a bijection from the quotient space
$\mathcal{D}^l_{\lambda,\mu}(M)/\mathcal{D}^{l-1}_{\lambda,\mu}(M)$ to
$\mathcal{S}^l_{\delta}(M)$.
Hence the space of symbols is nothing but
\[\mathcal{S}_{\delta}(M)
=\bigoplus_{l=0}^{\infty}\mathcal{S}^l_{\delta}(M),\]
endowed with the classical actions of $\mathrm{Diff}(M)$ and of $\mathrm{Vect}(M)$.
\subsection{Projective equivalence of connections}
We denote by $\mathcal{C}_M$ the space of torsion-free linear connections
on $M$. Two such connections are \emph{Projectively equivalent} if there
exists a one form $\alpha$ on $M$ such that their associated covariant
derivatives $\nabla$ and $\nabla'$ fulfill the relation
\[\nabla'_XY = \nabla_XY + \alpha(X) Y + \alpha(Y)X.\]
This equation was introduced by H. Weyl in \cite{Weyl}. He showed that it was a
 necessary and sufficient condition for the two connections to define the
same \emph{paths}, that is, the same geodesics up to parametrization. 
\subsection{Problem setting}
A \emph{quantization} on $M$ is a linear bijection $Q_M$ from
 the space of symbols
$\mathcal{S}_{\delta}(M)$ to the space of differential operators
$\mathcal{D}_{\lambda,\mu}(M)$ such that
\[\sigma (Q_M(S)) = S,\quad\forall S\in \mathcal{S}^k_{\delta}(M),\; \forall k
\in\mathbb{N}. \]
Roughly speaking, a \emph{natural quantization} is a quantization which
depends on a torsion-free connection and commutes with
 the action of diffeomorphisms. More precisely, a
natural quantization is a collection of maps (defined for every manifold $M$)
\[Q_M : \mathcal{C}_M \times \mathcal{S}_{\delta}(M)\to
\mathcal{D}_{\lambda,\mu}(M)\]
such that
\begin{itemize}
\item For all $\nabla$ in $\mathcal{C}_M$, $Q_M(\nabla)$ is a quantization,
\item If $\phi$ is a local diffeomorphism from $M$ to $N$,
  then one has
\[Q_M(\phi^*\nabla)(\phi^*S) = \phi^*(Q_N(\nabla)(S)),\quad\forall \nabla\in
  \mathcal{C}_N,\forall S\in \mathcal{S}_{\delta}(N).\]
\end{itemize}
A quantization $Q_M$ is \emph{projectively equivariant} if one has
$Q_M(\nabla) = Q_M(\nabla')$ whenever $\nabla$ and $\nabla'$ are projectively
equivalent torsion-free linear connections on $M$.

{\bf{Remark :}} The definition of a natural quantization was set by M. Bordemann
in functorial terms and  relates in this sense to the concept of natural
operators in differential geometry exposed in \cite{KMS}.
\section{Projective Cartan connections}
For the paper to be self-contained, we recall here the most important facts about Cartan connections. We begin
with a general definition and then we give more details about the projective
Cartan connections and their links with projective structures.
For more detailed information, the reader may refer to \cite{Kobabook}.
\subsection{Cartan connections}
Let $G$ be a Lie group and $H$ a closed subgroup. Denote by $\g$ and
$\h$ the corresponding Lie algebras. Let $P\to M$ be a principal
 $H$-bundle over $M$, such that $\mathrm{dim}\,M$ = $\mathrm{dim}\,G/H$.
 A Cartan
connection on $P$ is a $\g$-valued one-form $\omega$ on $P$ such that
\begin{itemize}
\item If $R_a$ denotes the right action of $a\in H$ on $P$, then
 $R_a^*\omega = Ad(a^{-1})\omega$,
\item If $k^*$ is the vertical vector
  field associated to $k\in\h$, then $\omega (k^*)= k$,
\item $\forall u\in P, \omega_u : T_uP\to \g$ is a linear bijection.
\end{itemize}
\subsection{Projective structures and Projective connections}
 We consider the group $G=PGL(m+1,\R)$ acting on $\R P^m$. We denote by $H$ 
the stabilizer of the element $[e_{m+1}]$ in
 $\R P^m$. One has 
\begin{equation}\label{hproj}H=\{\left(\begin{array}{cc}A & 0\\\xi & a
\end{array}\right) : A\in GL(m,\R),\xi\in \R^{m*}, a\not = 0\}/\R_0\mbox{Id},\end{equation}
and it follows that $H$ is the semi-direct product $G_0 \rtimes G_1$, where
$G_0$ is isomorphic to $GL(m,\R)$ and $G_1$ is isomorphic to $\R^{m*}$. 
Then there is a projection 
 \[\pi : H\to GL(m,\R) : \left[\left(\begin{array}{cc}A & 0\\\xi & a
\end{array}\right)\right]\mapsto \frac{A}{a}\]
The Lie algebra of  $PGL(m+1,\mathbb{R})$ is 
$gl(m+1,\mathbb{R})/\mathbb{R}Id$. It is thus  isomorphic to
$sl(m+1,\R)$ and it decomposes as a direct sum of subalgebras
 \[ \mathfrak{g}_{-1}\oplus\mathfrak{g}_{0}\oplus
\mathfrak{g}_{1}\cong\mathbb{R}^{m}\oplus gl(m,\mathbb{R})\oplus\mathbb{R}^{m*}.\]
The isomorphism is given by
\[\left[\left(\begin{array}{cc}A & v\\\xi & a
\end{array}\right)\right]\mapsto (v,A-a\,Id,\xi).\]
This correspondance induces a structure of Lie algebra on
$\mathbb{R}^{m}\oplus gl(m,\mathbb{R})\oplus\mathbb{R}^{m*}$. The Lie algebras
corresponding to $G_{0}$, $G_{1}$ and $H$ are respectively $\mathfrak{g}_{0}$,
$\mathfrak{g}_{1}$, and $\mathfrak{g}_{0}\oplus\mathfrak{g}_{1}$.

Let us denote by $G^2_m$ the group of 2-jets at the origin $0\in\R^m$ of local
 diffeomorphisms defined on a neighborhood of $0 $ and that
 leave $0$ fixed. The group $H$ acts on $\R^m$ by
linear fractional transformations that leave the origin fixed. This allows to
 view $H$  as a subgroup of $G^2_m$  , namely,
\begin{equation}\label{iso}
\i : H\to G^2_m : \left(\begin{array}{cc}A &0\\ \xi &
a\end{array}\right)\mapsto
(\frac{A^i_j}{a},-\frac{A^i_j\xi_k + A^i_k \xi_j}{a^2})
\end{equation}
A \emph{Projective
  structure on $M$} is then a reduction of the second order jet-bundle $P^2M$
to the group $H$. The following result (\cite[Prop 7.2 p.147]{Kobabook}) is the starting point of our method :
\begin{prop}[Kobayashi-Nagano]
There is a natural one to one correspondance between the projective
 equivalence classes of torsion-free linear connections on $M$
and the projective structures on $M$.
\end{prop}
In general, if $\omega$ is a Cartan connection defined on a $H$-principal bundle $P$, then its
curvature $\Omega$ is defined as usual by
 \begin{equation}\label{curv} 
\Omega = d\omega+\frac{1}{2}[\omega,\omega].
\end{equation}
The notion of \emph{Normal} Cartan connection is defined by natural conditions
imposed on the components of the curvature.

Now, the following result (\cite[p. 135]{Kobabook}) gives the relationship between  projective
structures and Cartan connections :
\begin{prop}
 A unique normal
 Cartan connection with values in the algebra $sl(m+1,\R)$ is
 associated to every projective structure $P$. This association is natural.
\end{prop}
The connection associated to a projective structure $P$ is called the normal
projective connection of the projective structure.
\subsection{Invariant differentiation}
We will use the concept of invariant differentiation with respect to
 a Cartan connection developed in \cite{Capinv,Capcart}. Let $P$ be a
projective structure and let $\omega$ be the associated normal projective
connection.
\begin{defi} Let $(V, \rho)$ be a representation of $H$. If $f\in
 C^{\infty}(P,V)$, then the invariant differential of $f$ with respect to
$\omega$ is the function \\
$\nabla^{\omega}f\in C^{\infty}(P,\R^{m*}\otimes V)$
defined by
\[ \nabla^{\omega}f(u)(X) = L_{\omega^{-1}(X)}f(u)\quad\forall u\in
P,\quad\forall X\in\R^m.\]
\end{defi}
We will also use an iterated and symmetrized version of the invariant 
differentiation 
\begin{defi}
If $f\in C^{\infty}(P,V)$ then $(\nabla^{\omega})^k f
\in C^{\infty}(P,S^k\R^{m*}\otimes V)$ is defined by
\[(\nabla^{\omega})^k f(u)(X_1,\ldots,X_k) = \frac{1}{k!}\sum_{\nu}
L_{\omega^{-1}(X_{\nu_1})}\circ\ldots\circ L_{\omega^{-1}(X_{\nu_k})}f(u)\]
for $X_1,\ldots,X_k\in\R^m$.
\end{defi}
\subsection{Lift of equivariant functions}
In order to make use of the invariant differentiation, we need to know the
relationship between equivariant functions on $P^1M$ and equivariant functions on
$P$. The following results were already quoted in \cite[p. 47]{Capinv}.

If $(V,\rho)$ is a representation of $\mathit{GL}(m,\R)$, then we define
a representation $(V,\rho')$ of $H$ by
\[\rho': H\to GL(V): \left[\left(\begin{array}{cc}A &0\\ \xi &
    a\end{array}\right)\right]\mapsto  \rho\circ \pi( \left[\left(\begin{array}{cc}A &0\\ \xi &
    a\end{array}\right)\right])
= \rho(\frac{A}{a})\]
for every $ A\in \mathit{GL}(m,\R), \xi\in
\R^{m*}, a\not=0$.

Now, using the representation $\rho'$,  we can give the relationship between
 equivariant functions on $P^1M$ and equivariant functions on $P$ :
If $P$ is  a projective structure on $M$, the natural projection $P^2M\to
P^1M$ induces a projection $p :P\to P^1M$ and we have :
\begin{prop} If $(V,\rho)$ is a representation of $GL(m,\R)$, then the map
\[p^* : C^{\infty}(P^1M,V)\to C^{\infty}(P,V) : f \mapsto f\circ p\]
defines a bijection from $C^{\infty}(P^1M,V)_{\mathrm{GL}(m,\R)}$ to
$C^{\infty}(P,V)_{H}$.
\end{prop}
This result is well-known and comes from the following  facts 
\begin{itemize}
\item $(p,Id,\pi)$ is a
morphism of principal bundles from $P$ to $P^1M$ 
\item the equivariant
functions on $P$ are constant on the orbits of the action of $G_1$ on $P$.
\end{itemize}
Now, since $\R^m$ and $\R^{m*}$ are natural representations of $GL(m,\R)$,
they become representations of $H$ and we can state an important property of the invariant differentiation :
\begin{prop}\label{gonabla}If $f$ belongs to $C^{\infty}(P,V)_{G_0}$ then
  $\nabla^{\omega}f\in C^{\infty}(P,\R^{m*}\otimes V)_{G_0}$.
\end{prop}
\begin{proof}The result is a direct consequence of the Ad-invariance of the
  Cartan connection $\omega$. 
\end{proof}
The main point that we will discuss in the next sections is that this result
is not true in general for $H$-equivariant functions : for an $H$-equivariant
function $f$, the function
$\nabla^{\omega}f$ is in general not $G_1$-equivariant.

As we continue, we will use the representation $\rho'_*$ of the Lie algebra of
$H$ on $V$. If we recall that this algebra is isomorphic to
 $gl(m,\R)\oplus \R^{m*}$
then we have 
\begin{equation}\label{rho}\rho'_* (A, \xi) = \rho_*(A),\quad\forall A\in
  gl(m,\R), \xi\in \R^{m*}.\end{equation}
In our computations, we will make use of the infinitesimal version of the
equivariance relation : If $f\in C^{\infty}(P,V)_H$ then one has
\begin{equation}\label{Invalg}
L_{h^*}f(u) + \rho'_*(h)f(u)=0,\quad\forall h\in gl(m,\R)\oplus\R^{m*}\subset sl(m+1,\R),
\forall u\in P.
\end{equation}
\section{Construction of a projectively invariant quantization}
The existence of a natural and projectively equivariant quantization is linked
to the existence of an $sl(m+1,\R)$-equivariant quantization in the sense of
\cite{LO} in the flat situation. It is known
that for some \emph{critical values} of $\delta$ such a quantization does
 not exist. Let us recall the
 following definition of \cite[Prop 2, p. 289]{Lecras} :
\begin{defi}
We define the numbers 
\[\gamma_{2k-l} = \frac{(m+2 k - l -(m+1)\delta)}{m+1}.\]
A value of
$\delta$ is critical if there exists $k,l\in \N$ such that $1\leq l\leq k$ and $\gamma_{2k-l}=0$.
\end{defi}
One of the results of \cite{Lecras} is then the following
\begin{thm}\label{plat1}
If $\delta$ is not critical, there exist an unique $sl(m+1,\R)$-equivariant
quantization.
\end{thm}
A link between the natural quantization and the $sl(m+1,\R)$-equivariant
 quantization is then given in \cite{Leconj} by 
\begin{thm}\label{plat2}
If $Q_{M}$ is a natural projectively equivariant quantization and if we denote
by $\nabla_0$ the flat connection on $\R^m$, then
$Q_{\mathbb{R}^{m}}(\nabla_{0})$ is $sl(m+1,\R)$-equivariant.
\end{thm}
Now, let us introduce the divergence operator associated to a Cartan
connection. This operator will be the main tool of our construction.
\subsection{The Divergence operator}
We fix a basis $(e_1,\ldots,e_m)$ of $\R^m$ and we denote by
$(\epsilon^1,\ldots,\epsilon^m)$ the dual basis in $\R^{m*}$.

The \emph{Divergence operator} with respect to the Cartan connection
$\omega$ is then defined by
\[div^{\omega} : C^{\infty}(P,S^k_{\delta}(\R^m))
\to C^{\infty}(P,S^{k-1}_{\delta}(\R^m)) :
S\mapsto \sum_{j=1}^m i(\epsilon^j)\nabla^{\omega}_{e_j}S,\]
where $i$ denotes the inner product.

This operator can be seen as a curved generalization of the
divergence operator
used in \cite{LO}.
The following propositions shows its most important properties.
\begin{lem}\label{goinv}
If $S\in C^{\infty}(P,S^k_{\delta}(\R^m))_{G_0}$
 then $div^{\omega}S
\in C^{\infty}(P,S^{k-1}_{\delta}(\R^m))_{G_0}$.
\end{lem}
\begin{proof} 
This can be checked directly from the definition.
One can also remark that $div^{\omega}S$ is the contraction of the
invariant function $\nabla^{\omega}S$ (see proposition~\ref{gonabla}) and
 of the constant and invariant function
\[ID : P\to \R^m\otimes \R^{m*},
u\mapsto\sum_{j=1}^m \epsilon^j\otimes e_j.\]
\end{proof}
The purpose of the next results is to measure the failure of $G_1$-
equivariance of the operators defined so far. At the infinitesimal level, in
view of equations (\ref{rho}) and (\ref{Invalg}), this
leads to the computation of the commutator of these operators with the Lie
derivative $L_{h^*}$, for $h\in\g_1$.
 We begin with the divergence operator :
\begin{lem}\label{div1}
For every $S\in C^{\infty}(P,S^k_{\delta}(\R^m))_{G_0}$
we have
\[L_{h^*} div^{\omega}S - div^{\omega}L_{h^*}S = (m+1) \gamma_{2k-1} i(h)S,\]
for every $h\in \R^{m*}\cong \g_1$.
\end{lem}
\begin{proof}
First we remark that the Lie derivative with respect to a vector field
commutes with the evaluation :
If $\eta^1,\ldots,\eta^{k-1}\in \R^{m*}$, we have
\[\begin{array}{lll}(L_{h^*}div^{\omega}S)(\eta^1,\ldots,\eta^{k-1})
& =& L_{h^*}(div^{\omega}S(\eta^1,\ldots,\eta^{k-1}))\\
 & = &
\sum_{j=1}^m(L_{h^*}L_{\omega^{-1}(e_j)}S(\epsilon^j,\eta^1,\ldots,\eta^{k-1}))\end{array}\]
Now, the definition of a projective Cartan connection implies 
the relation
 \[[h^*,\omega^{-1}(X)] = \omega^{-1}([h,X]),\quad\forall h\in
 gl(m,\R)\oplus \R^{m*},X\in\R^m,\]
where the bracket on the right is the one of $sl(m+1,\R)$.
It follows that the expression we have to compute is equal to
\[\sum_{j=1}^m( L_{\omega^{-1}(e_j)}L_{h^*}S(\epsilon^j,\eta^1,\ldots,\eta^{k-1})
+ (L_{[h,e_j]^*}S)(\epsilon^j,\eta^1,\ldots,\eta^{k-1})).\]
Finally, using the relation (\ref{Invalg}), we obtain
\[\begin{array}{lll}&& div^{\omega}(L_{h^*}S) -(\rho_{*}'([h,e_j])S)(\epsilon^j,\eta^1,
\ldots,\eta^{k-1})\\
&=&div^{\omega}(L_{h^*}S) -(\rho_{*}(h\otimes e_j + \langle h,e_j\rangle Id)S)(\epsilon^j,\eta^1,
\ldots,\eta^{k-1}).\end{array}\]
The result then easily follows from the definition of $\rho$ on $S^k_{\delta}(\R^m)$.
\end{proof}
Eventually we obtain
\begin{prop}\label{div2}
For every $S\in C^{\infty}(P,S^k_{\delta}(\R^m))_{G_0}$,
we have
\[L_{h^*} (div^{\omega})^lS - (div^{\omega})^lL_{h^*}S = (m+1)l \gamma_{2k-l}
i(h) (div^{\omega})^{l-1}S,\]
for every $h\in \R^{m*}\cong \g_1$.
\end{prop}
\begin{proof}
For $l=1$, this is just lemma \ref{div1}. Then the result follows by
induction,
 using lemmas \ref{goinv} and \ref{div1}.
\end{proof}
Next, we analyse the failure of invariance of the iterated invariant differentiation :
\begin{prop}\label{nablag1}
If $f\in\mathcal{C}^{\infty}(P,\Delta^{\lambda}\mathbb{R}^{m})_{G_0}$, then 
\[L_{h^{*}}(\nabla^{\omega})^{k}f - (\nabla^{\omega})^{k}L_{h^{*}}f
 =-k((m+1)\lambda+k-1)((\nabla^{\omega})^{k-1}f\vee
h),\]
for every $h\in \R^{m*}\cong \g_1$.
\end{prop}
\begin{proof}If $k=0$, then the formula is obviously true. Then we proceed
 by induction. In view of the symmetry of the expressions that we have to
 compare, it is sufficient to check that they coincide when evaluated on the
$k$-tuple $(X,\ldots,X)$ for every $X\in\R^m$. The proof is similar to the one
 of lemma \ref{div1} : first the evaluation and the Lie derivative commute :
\[(L_{h^{*}}(\nabla^{\omega})^{k}f)(X,\ldots,X)
= L_{h^*}((\nabla^{\omega})^{k}f)(X,\ldots,X)).\]
Next, we use the definition of the iterated invariant differential and we let
the operators $L_{h^*}$ and $L_{\omega^{-1}(X)}$ commute so that the latter
expression becomes
\[L_{\omega^{-1}(X)}L_{h^{*}}((\nabla^{\omega})^{k-1}f)(X,\ldots,X)
+(L_{[h,X]^{*}}((\nabla^{\omega})^{k-1}f))(X,\ldots,X).\]
By the induction, the first term is equal to 
\[(\nabla^{\omega})^{k}L_{h^{*}}f(X,\ldots,X)  -(k-1)((m+1)\lambda+k-2)((\nabla^{\omega})^{k-1}f\vee
h)(X,\ldots,X).\]
For the second term, we use proposition \ref{gonabla} and relation
(\ref{Invalg}) and we obtain
\[(\rho_{*}((h\otimes X) +\langle h, X\rangle Id)((\nabla^{\omega})^{k-1}f))(X,\ldots,X).\]
The result follows by the definition of $\rho_*$.
\end{proof}
\subsection{The main result}
In this section, we give an explicit formula for the natural and projectively
equivariant quantization, using the properties of the iterated invariant
differentiation and of the divergence operator.
\begin{thm} If $\delta$ is not critical, then the
 collection of maps \\
$Q_M : \mathcal{C}_M\times \mathcal{S}_{\delta}(M)\to
\mathcal{D}_{\lambda,\mu}(M)$ defined by
\begin{equation}\label{formula}Q_M(\nabla, S)(f) = p^{*^{-1}}(\sum_{l=0}^k C_{k,l} \langle Div^{\omega^l}
p^*S,\nabla_s^{\omega^{k-l}}p^*f\rangle),\forall S\in \mathcal{S}^k_{\delta}(M)\end{equation}
defines a projectively invariant natural quantization if 
\[C_{k,l} =\frac{(\lambda + \frac{k-1}{m+1})\cdots (\lambda +
  \frac{k-l}{m+1})}{\gamma_{2k-1}\cdots
  \gamma_{2k-l}}\left(\begin{array}{c}k\\l\end{array}\right),\forall l\geq 1,\quad C_{k,0}=1.\]
\end{thm}
\begin{proof}
First, we have to check that the formula makes sense : the function
\begin{equation}\label{Hinv}\sum_{l=0}^k C_{k,l} \langle (Div^{\omega})^l
p^*S,(\nabla^{\omega})^{k-l}p^*f\rangle\end{equation}
has to be $H$-equivariant.
It is obviously $G_0$-equivariant by proposition \ref{gonabla} and lemma
\ref{goinv}. It is then sufficient to check that it is $\g_1$-equivariant. This follows directly from propositions
\ref{div2} and \ref{nablag1}  and from the relation 
\begin{equation}\label{ckl}C_{k,l}l(m+2k-l-(m+1)\delta)=C_{k,l-1}(k-l+1)((m+1)\lambda+k-l).\end{equation}
Next we see, using the results of \cite[p.47]{Capinv}
 that the principal symbol 
of $Q_M(\nabla, S)$ is exactly $S$, and formula (\ref{formula}) defines a
quantization, that is projectively invariant,
 by the definition of $\omega$.
Next, the naturality of the quantization defined in this way is easy to
understand : it follows from the naturality of the association of a
projective structure $P\to M$ endowed with a normal Cartan connection
 $\omega$ to a class of projectively equivalent torsion-free connections on $M$
  and from the naturality of the lift of the equivariant
 functions on $P^1M$ to equivariant functions on $P$.
\end{proof}
\noindent{\bf{Remarks :}}
\begin{itemize}
\item Theorems \ref{plat1} and \ref{plat2} directly imply that, when
$M$ is taken to be $\R^m$ and $\nabla$ is the flat connection, formula
\ref{formula} must coincide with the ones of \cite{LO} (formulas 4.14 and
4.15) and \cite{DO} (Formula 2.4), at least when $\delta$ is not critical.
 What is more surprising is that our
coefficients $C_{k,l}$ coincide with the ones of \cite{DO} (formulas 2.5 and
3.6), up to a combinatorial coefficient, which is due to a slightly different
definition of the divergence operator. In particular, our formula can be
expressed, as the one of \cite{DO}, in terms of hypergeometric functions.
\item A long computation, involving the explicit form of the normal Cartan
  connection in coordinates allows to show that our formula coincides for the
  case of third order differential operators with the formula provided by
 S. Bouarroudj in \cite{Bou1,Bou2}.
\item It was shown in \cite{LO} how an $sl(m+1)$-equivariant quantization
  (in the flat case) induced a projectively invariant star product on the space
  of symbols. It was shown in \cite{Brylinski} that the bilinear
  operators appearing in the deformation were not bidifferential
  operators. We conjecture that our formula \ref{formula} 
will induce a deformation of the
  algebra of symbols (depending on a connection) and that this deformation will not be local and therefore
  will not be a star-product in the sense of \cite{BFFLS}.

\end{itemize}
Now, the proof of the previous theorem also allows to analyse the existence
problem when $\delta$ is a critical value :
assume that there exist $k\in \N$ and $r\in \N$ such that $1\leq r \leq k$ and
$\gamma_{2k-r} = 0$. Then if there exists $i\in\{1,\ldots,r\}$ such that
 $\lambda=-\frac {k-i}{m+1},$ then one can replace the coefficients
$C_{k,i},\ldots,C_{k,k}$ by zero and the function (\ref{Hinv}) is still
$H$-equivariant. Then the collection $Q_M$ still defines a
projectively equivariant and natural quantization. If $\lambda$ does not
belong to the set $ \{-\frac {k-1}{m+1},\ldots,-\frac {k-r}{m+1}\}$, then
there is no solution since the $sl(m+1,\R)$- equivariant quantization in the
sense of \cite{LO,Lecras} does not exist. 
To sum up, we have shown the following
\begin{thm}
There exists a natural and projectively equivariant quantization if and only
if there exists an $sl(m+1,\R)$-equivariant quantization in the sense of
\cite{LO} over $M=\R^m$.
\end{thm}
\section*{Acknowledgements}
It is a pleasure to thank S. Hansoul and P. Lecomte for numerous fruitful
discussions and for their interest in our work. We thank D. Sternheimer for
his comments and suggestions.
We are grateful to C. Roberts, D.J. Saunders, M. Crampin and W. Bertram 
for their help and the references they gave us about Thomas-Whitehead and
Cartan connections.

F. Radoux thanks the Belgian FRIA for his Research Fellowship.

\bibliographystyle{plain} \bibliography{prquant}

\end{document}